\documentclass[12pt]{article}
\usepackage{amsmath, amsfonts, amssymb, amsthm} \parskip1pt
\newtheorem{theorem}{Theorem}[section]
\newtheorem{proposition}[theorem]{Proposition}

\newtheorem{remark}[theorem]{Remark}

\newtheorem{example}[theorem]{Example}

\newtheorem{corollary}[theorem]{Corollary}

\def\R{\mathbb R} \def\Z{\mathbb Z}

\def\N{\mathbb N}

\def\norm{|\!|}

\def\Ad{\hbox{\rm Ad}}

\def\norm{|\! | }

\begin{document}

\iffalse
\documentclass[12pt]{amsart}
\usepackage{graphicx}
%\usepackage{latexsym}
%\usepackage[all]{xy}
%\usepackage{vpage}
\usepackage{graphics, setspace}
\usepackage{amsmath,amsthm, amsfonts,amssymb, stmaryrd,yfonts,pxfonts,pifont,eufrak, bbm}
%\usepackage{epsfig,enumerate}
%\usepackage{indentfirst}
%\usepackage{makeidx}
% ----------------------------------------------------------------
\vfuzz2pt % Don't report over-full v-boxes if over-edge is small
\hfuzz2pt % Don't report over-full h-boxes if over-edge is small
% THEOREMS -------------------------------------------------------

\fi

%\newtheorem{th}{Theorem}
\newtheorem{Theorem}{Theorem}[section]
\newtheorem{Obs}[Theorem]{Observation}
\newtheorem{Cor}{Corollary}
\newtheorem{ex}{Example}
\newtheorem{p}{Proposition}
\theoremstyle{definition}
\newtheorem{defn}{Definition}
\theoremstyle{remark}
\newtheorem{Remark}{Remark}
\numberwithin{equation}{subsection}

\iffalse
\newcommand{\norm}[1]{\left\Vert#1\right\Vert}
\newcommand{\abs}[1]{\left\vert#1\right\vert}
\newcommand{\set}[1]{\left\{#1\right\}}
\newcommand{\R}{\mathbb R}
\newcommand{\eps}{\varepsilon}
\newcommand{\ra}{\rightarrow}
\newcommand{\To}{\longrightarrow}
\newcommand{\BX}{\mathbf{B}(X)}
\newcommand{\A}{\mathcal{A}}
\newcommand{\N}{\mathbb N}
\newcommand{\Z}{\mathbb Z}
\fi
\newcommand{\T}{\mathbb T}
\newcommand{\U}{\mathcal U}

\let\Ga=\Gamma
\let\ap=\alpha

\title{Affine Almost Automorphic Actions on Compact Nilmanifolds}

\author{S.G.~Dani, Riddhi Shah \& Puneet Sharma}

%\begin{document}

\maketitle
\date{}

\begin{abstract}
We discuss conditions under which an  affine automorphism of a compact
nilmanifold is almost automorphic, and the structure of such automorphisms
from a dynamical point of view.

\end{abstract}

\section{Introduction}

Let $X$ be a compact second countable space and $\alpha$ a homeomorphism
of $X$.
A point $x\in X$ is said to be  {\it almost automorphic} for the
action of $\alpha$ if for any $x, y \in X$ and any
sequence $\{k_i\}$ in $\Z$ such that $\alpha^{k_i}(x)\to y$ we
have $\alpha^{-k_i}(y)\to x$.
The action of $\alpha$ on $X$ is said to be {\it almost automorphic}
if every $x\in X$ is almost automorphic for the $\alpha$-action. It may
be expected that such a behaviour would be
rare. In this note we discuss the notion  for affine automorphisms
of compact nilmanifolds; the results
illustrate the restrictive nature of the condition.

Let $N$ be a simply connected nilpotent (Lie) group and $\Ga$
be a lattice in $N$ (i.e.\ $\Ga$ is a discrete co-compact subgroup in $N$).
Then $N/\Gamma$ is a compact manifold, and the quotients arising in
this way are called nilmanifolds. By an automorphism $A$ of $N/\Ga$ we
mean a self-map such that $A(g\Ga)=\tilde A(g)\Ga$ for all $g\in N$,
where $\tilde A$ is a (Lie) automorphism of
$N$ leaving $\Ga$ invariant (viz. such that $\tilde A(\Ga)=\Ga$). For
$a\in N$, the
map $T_a:N/\Ga \to N/\Ga$ given by $g\Ga \mapsto ag\Ga$, for all $g\in
G$, is called the
translation by $a$. A transformation of $N/\Ga$ of
the form $T=T_a\circ A$ where $a\in N$ and $A$ is an
automorphism of $N/\Ga$ is called {\it affine automorphism} of $N/\Ga$.
An automorphism $\tilde A$ of $N$ is said to be unipotent if the
corresponding automorphism of the Lie algebra is a
unipotent linear transformation, and we call an automorphism $A$
of $N/\Ga$ unipotent if it is induced by a unipotent automorphism
$\tilde U$ such that $\tilde U(\Ga)=\Ga$.

Our aim is to characterize affine automorphisms of compact
nilmanifolds which are almost automorphic. It turns out that
if $T=T_a\circ A$ is almost automorphic then for some $r\in\N$, $A^r$ is
unipotent (see Proposition~\ref{general}, and the remark preceding
it). We note also that for a transformation $T$,  $T^r$ is almost
automorphic, for $r\in \N$, if and only if $T$ is almost automorphic.
Since $T^r$ is of the form $T_{a_r}\circ A^r$, for some $a_r\in N$,
this shows that it
suffices to consider affine automorphisms
$T=T_a\circ U$, with $U$ unipotent, and in the interest of simplicity
of exposition, until \S 6  we shall restrict to this situation.

We prove the following.

\begin{theorem}\label{thm:main}
Let $N$ be a simply connected nilpotent Lie group and $\Gamma$
be a lattice in $N$.   Let $T=T_a\circ  U$ be an affine
automorphism of  $N/\Ga$, where $U$ is a unipotent automorphism of
$N/\Ga$. Let  $\tilde U$ be the automorphism of $N$ corresponding
to $U$. The point $\Ga$ is almost automorphic for
the $T$-action if and only if there exists a closed abelian subgroup $A$
of $N$
such that $A\Ga$ is closed, $a\in A$ and $\tilde U(g)=g$ for all $g
\in A$; in particular, if $\Ga$ is almost automorphic then $\tilde U(a)=a$.
\end{theorem}

The theorem implies the following towards unipotent affine
automorphisms being almost automorphic.

\begin{corollary}\label{cor:main}
Let the notation be as in Theorem~\ref{thm:main}.
Then $T$ is almost automorphic if and only if there exists a closed
connected abelian subgroup $A$ of $N$ such that $A\Ga$ is closed,
$x^{-1}a\tilde U(x)\in A$ for all $x\in N$ and $\tilde U(g)=g$ for all
$g\in A$. In particular, if $T$ is minimal and almost automorphic then
$N$ is abelian and $T$ is a translation of the torus $N/\Ga$.
\end{corollary}

The ``if'' part is straightforward in both Theorem~\ref{thm:main} and
Corollary~\ref{cor:main} (see \S \ref{mainthm} for
details). The main point of the results is that the converse also holds.

The conclusion of the corollary means that, from the point of view of
dynamics, almost automorphic affine transformations
of compact nilmanifolds consist of a collection of  translations of
tori put together; the ``putting
together'' would of course depend on the nilmanifold and the
transformation.

For translations we get the following characterisation.

\begin{corollary}\label{transl}
Let $N$ be a simply connected nilpotent Lie group and $\Gamma$ be a
lattice in $N$. Let $a\in N$. Then the translation $T_a$ of $N/\Ga$ is
almost automorphic if and only if there exists a closed connected
abelian normal subgroup $A$
of $N$ such that $A\Gamma$ is closed.

\end{corollary}

The following describes certain simple necessary conditions, at the
``infinitesimal'' (or Lie algebra) level, on the
automorphism and the translating element for an affine automorphism of
a compact nilmanifold to be almost automorphic.

\begin{corollary}\label{cor:Lie}
Let the notation be as in Theorem~\ref{thm:main}.  Furthermore let
$\mathfrak N$ be the Lie algebra of $N$, $\mathfrak
U$ be the Lie automorphism of $\mathfrak N$ corresponding to $\tilde U$,
and $\Ad_a$ be the adjoint transformation of $\mathfrak N$ corresponding
to $a$. Let $I$ and $0$ denote the identity and zero transformations
of $\mathfrak N$, respectively.  Suppose
$T$ is almost automorphic. Then the following holds:

i) $(\mathfrak U-I)(\Ad_a\circ \mathfrak U-I)=0$; in particular, if $T$ is an
automorphism (viz. $a$ is the identity) then $\tilde U$ has nilrank at
most 1, namely $(\tilde U-I)^2=0$.

ii) The image of $\Ad_a\circ \mathfrak U-I$ is an abelian Lie subalgebra.
\end{corollary}

As we shall see from an example these conditions are not sufficient
for the converse to hold, even in the case of automorphisms.

\section{Minimal Translations}

In this section we shall prove the following proposition which
in particular proves the statement of Theorem~\ref{thm:main} in the
special case of minimal translations.  We recall that
an action is said to be {\it minimal} if
there is no proper nonempty closed invariant subset, or
equivalently if every orbit is dense.

\begin{proposition}\label{minimal}
Let $N$ be a simply connected nilpotent Lie group
  and $\Ga$ be a lattice in $N$. Let $\Phi=\{a_t\}$ be a one-parameter
  subgroup of $N$. Suppose that the action of $\Phi$ on $N/\Ga$ is
  minimal, and that $\Ga$ is an almost automorphic point of the
  action.  Then $N$ is abelian.
\end{proposition}

\proof We proceed by induction on the dimension $n$ of $N$. The
assertion is evident for $n=1$. Now suppose the
dimension is $n$ and that the proposed statement holds when the
dimension is less than $n$.

Let $Z$ denote the centre of $N$. Then $Z$ is topologically isomorphic
to $\R^k$ for some $k$ and $Z\cap \Ga$ is a lattice in $Z$. Let
$P$ be a one-parameter subgroup of $Z$ containing a nontrivial element
of $\Ga$; we note that $P\Ga$ is a subgroup, and as $P/(P\cap \Ga)$ is
compact, $P\Ga$ is closed. Now,
$N/P$ is a simply connected nilpotent Lie group and $\Ga P/P$ is a
lattice in $N/P$.
The $\Phi$-action on $N/\Ga$ factors to an action on $(N/P)/(\Ga P/P)$
which is minimal and almost automorphic. By the induction hypothesis
therefore $N/P$ is abelian.

Now suppose, if possible, that $N$ is nonabelian. Let $\mathfrak N$ be
the Lie algebra of $N$ and $\exp :\mathfrak N \to N$ be the exponential
map. Let $\Lambda$ be the (additive) subgroup of $\mathfrak N$ generated
by $ \exp^{-1}(\Ga)$. Then $\Lambda$ is a lattice in $\mathfrak N$ and
$\exp (\Lambda)$ is a lattice in $N$; (see \cite{raghunathan},
Theorem~2.12, addendum to the main statement). Clearly $\exp
(\Lambda) $ contains $\Ga$ as a subgroup of finite index. Moreover,
since $\Gamma$ is $\tilde U$-invariant, so is $\exp (\Lambda)$, and
the factor of $T$ on $N/\exp (\Lambda) $ is almost automorphic and
minimal. Therefore, in the proof of the theorem replacing $\Ga$ by
$\exp (\Lambda)$ and modifying notation we may assume that $\Ga
=\exp (\Lambda)$, with $\Lambda$ a lattice in $\mathfrak N$. Now let
$\mathfrak P$
be the (one-dimensional) Lie subalgebra corresponding to $P$. We
note that $\mathfrak P \cap \Lambda$ is an infinite cyclic subgroup and
choose a generator $e_0$. Then $e_0$ can be extended to a basis of
$\Lambda$, as a free abelian group, say $\{e_0,e_1, \dots, e_d\}$.
Let $W$ be the subspace of $\mathfrak N$ spanned by $\{e_1, \dots,
e_d\}$. Since $N/P$ is abelian, $\mathfrak N$ can be realised as $W+\R$
with the Lie bracket operation given by
$[(w_1,t_1),(w_2,t_2)]=[0,\alpha (w_1,w_2)]$ for all $w_1,w_2 \in W$
and $t_1,t_2\in \R$, where $\alpha$ is a bilinear form on
$W$. Correspondingly $N$ can realised as $W\times
\R$, with the multiplication given by $(w_1,t_1)\cdot
(w_2,t_2)=(w_1+w_2, t_1+t_2+\mu (w_1,w_2))$, for all $w_1,w_2\in W$,
$t_1,t_2\in \R$, with $\mu$ a bilinear form on $W$, and the
exponential map $\exp : \mathfrak N \to N$ is given by $\exp ((w,t))=(w,
\frac 12 \mu (w,w) +t)$ for all $w\in W$ and $t\in \R$.
%We note that $\mu(\Delta \times \Delta)$ is contained in a discrete
%subgroup, and modifying the scale on $\R$ if necessary, we assume that
%consists of integers.

Let $\Delta$ be the integral lattice in $W$ with respect to the basis
$\{e_1, \dots ,e_d\}$. Then $\Delta$ is the image of $\Ga$ when $W$ is
viewed as $N/P$, and we note also that for any $\delta \in \Delta$, as
a subset of $\mathfrak N$, we have
$\exp (\delta)=(\delta ,\frac 12 \mu (\delta, \delta)) \in \Ga$.

Let $$
F=\{u\in W\mid u=\sum_{i=1}^d\alpha_i e_i, \hbox{ \rm with }
|\alpha_i| \leq \frac12 \hbox{ \rm for all } i\}.$$
Consider any $w\in W$. Then there exist $u\in F$ and  $\delta \in
\Delta $ such that $w=u+\delta$.  Then in $N$, we have $(w,0)=
(u,-\mu (u,\delta)-\frac 12 \mu(\delta, \delta))\cdot (\delta,
\frac 12 \mu (\delta, \delta))$, and since $(\delta,
\frac 12 \mu(\delta, \delta))\in \Ga$ we get that $(w,0)\Ga=
(u,-\mu (u,\delta)-\frac 12 \mu(\delta, \delta))\Ga$.
For $w\in W$ such that $u$ is in the interior of $F$ the $\delta$ as
above is unique and
in this case we shall denote the element $(-\mu (u,\delta)-\frac 12
\mu(\delta, \delta))$ by $\rho (w)$; hence $(w,0)\Ga = (u,
\rho (w))\Ga$.
We note also that given such a $w$, we have $(-w,0)=(-u,-\mu
(-u,-\delta)-\frac 12 \mu(-\delta, -\delta)))\cdot (-\delta,
\frac 12 \mu (-\delta,- \delta)))$, so $\rho (-w)= -\mu
(-u,-\delta)-\frac 12 \mu (-\delta,-\delta)))=\rho (w)$.
%(Here, $(-w,0)\Ga=(-u,\rho(-w))\Ga=(-u,\rho(w))\Ga$.)

Since the $\Phi$-action is minimal it follows that there exists $r\in
\R$ such that $T_{a_r}$ is minimal; this may be seen from Parry's
theorem \cite{parry} together with the fact that it is true for tori.
Let $(b,s)$, with $b\in W$ and $s\in \R$ be the pair corresponding to
$a_r\in  \Phi \subset N$, and let $T=T_{a_r}$. Then for any $w\in W$
and $t\in \R$ we have  $T((w,t)\Ga)=(b+w, s+t + \mu (b,w))\Ga.$
Recursively we see that
$$T^k((w,t)\Ga) = (kb+w, ks+t+\sum_{j=0}^{k-1} \mu (b,jb+w)) \Ga$$ for
all $k\geq 1.$ Also, for all $w\in W$ and $t\in \R$ we have
$T^{-1}((w,t)\Ga)= (-b+w, -s+t-\mu(b,-b+w))\Ga$ and recursively,
 $$ T^{-k}((w,t)\Ga) = (-kb+w, -ks+t-\sum_{j=1}^{k} \mu (b,-jb+w)) \Ga$$
for all $k\geq 1.$

Let $\xi\in \R$ be irrational.
We now choose a sequence $\{k_i\}$ in $\N$ such that $T^{k_i}(\Ga) \to
(0,\xi)\Ga$, as follows. Let $S$ be the
translation $T_{(b,0)}$ of $N/\Ga$ and $\tau$ be the (cartesian)
product transformation $S\times T_{\mu(b,b)}$  of
$(N/\Ga)\times (\R/\Z)$; viz $\tau$  consists of translations by
$(b,0)$ and $\mu
(b,b)$ in the respective components. If  $\tau$ is
minimal, then we can choose $\{k_i\}$ such that $S^{k_i}(\Ga)\to
(0,\xi)\Ga$ and $k_i\mu(b,b)\to 0$ mod~$\Z$. Now suppose $\tau$ is not
minimal. We note that since the factors of $S$ and $T$ on $W/\Delta$
coincide, by Parry's theorem \cite{parry} $S$ is minimal. In this case
we choose a $\{k_i\}$ such that  $S^{k_i}(\Ga)\to
(0,\xi)\Ga$. In either case, passing to a subsequence
we shall also assume that $(0,k_is)\Ga\to (0,\sigma)$ as $i\to
\infty$, where $\sigma \in \R$.

We note the following in the case when $\tau$ is
not minimal. Under this condition,
$T_b\times  T_{\mu(b,b)}$ of $(W/\Delta) \times (\R/\Z)$ is also not
minimal. This further implies that there exist $p \Z$, $q\in \N$, and
a linear form $\psi$ on
$W$ such that $\psi (\Delta)\subset \Z$ and $\mu (b,b)=(\psi
(b)+p)/q$. Since $S^{k_i}(\Ga)\to (0,\xi)\Ga$, $\{T_b^{k_i}(0)+\Delta\}$
which is the same as $\{k_ib+\Delta\}$ converges to the identity
element in $W/\Delta$. Since $\mu (b,b)=(\psi
(b)+p)/q$, with $\psi$ and $p,q$ as above, it follows that every limit
point of $\{k_i\mu (b,b)+\Z\}$ in $\R/Z$ is a $q$th root of
unity, and passing to a subsequence we may assume that $\{k_i\mu
(b,b)+\Z\}$ converges to a root of unity.

Thus in either of the cases, we may assume $k_i\mu (b,b)+\Z \to \eta
+\Z$, with $\eta$ rational.

Now, for all $i$ we have, $T^{k_i}(\Ga)=(k_ib,
k_is+\theta_i)\Ga$, and $S^{k_i}(\Ga)=(k_ib,\theta_i)\Ga$, where $\theta_i=\sum_{j=0}^{k_i-1} \mu
(b,jb))=\frac12 (k_i-1)k_i\mu (b,b)$. Let  $u_i\in F$ and $\rho_i\in \R$
be such that $(k_ib,0)\Ga=(u_i, \rho_i)\Ga$. Then for all $i$,
$T^{k_i}(\Ga) = (u_i,k_is+\rho_i+\theta_i)\Ga$  and $S^{k_i}(\Ga) =
(u_i,\rho_i+\theta_i)\Ga$, and since the latter sequence converges to
$(0,\xi)\Ga$ we get that $u_i\to 0$ and  $\rho_i+ \theta_i\to
\xi$  mod $\Z$, as $i\to \infty$. Since $(0,k_is)\Ga\to (0,\sigma)$ as
$i\to \infty$ this shows that  $T^{k_i}(\Ga)\to
(0,\sigma+\xi)\Ga$, as $i\to \infty$.

On the other hand, for all $i$,  $T^{-k_i}(\Ga)=(-k_ib,
-k_is+\theta_i+k_i\mu(b,b))\Ga$. Since $u_i\to 0$, for all large
$i$ we have $\rho_i=\rho (k_ib)$, and since $\rho (k_ib)=\rho (-k_ib)$
we get that
$(-k_ib,0)\Ga=(-u_i, \rho_i)\Ga$ for all large $i$. Hence
$T^{-k_i}(\Ga)=(-u_i, -k_is+\rho_i+\theta_i+k_i\mu(b,b))\Ga$, for all
$i$. Using the convergences as above we see from this that as $i\to
\infty$, $T^{-k_i}(\Ga)\to (0,-\sigma+ \xi+\eta)\Ga$. Hence $T^{-k_i}
((0,\sigma+\xi)\Ga)\to (0, 2\xi+\eta)\Ga$.  Since $\xi$ is
irrational and $\eta$ is rational  $(0, 2\xi+\eta)\Ga \neq \Ga$; this
is a contradiction to the
condition that $\Ga$ is an almost automorphic point for the
$\Phi$-action.
This shows that $N$ must be abelian, thus proving the
proposition. \hfill $\Box$

\section{Proof of Theorem~\ref{thm:main} for translations}

We next prove the Theorem~\ref{thm:main} for all
translations. Specifically we prove the following.

\begin{proposition}\label{prop}
Let $N$ be a simply connected nilpotent Lie group and $\Ga$ be a
lattice in $N$. Let $\Phi= \{a_t\}$ be a one-parameter
subgroup of $N$ and let $a=a_t$ for some $t\neq 0$. Then the following
conditions are equivalent.

i) $\Ga$ is almost automorphic for the action of $T_a$ on $N/\Ga$;

ii) $\Ga$ is almost automorphic for the $\Phi$-action on $N/\Ga$;

iii)  there exists a closed connected abelian subgroup
$A$ of $N$ such that $\Phi$ is contained in $A$ and $A\Ga$ is closed.
\end{proposition}

\proof The proof of equivalence of (i) and (ii) is routine and we
shall omit it.

ii) $\implies$ iii): We  proceed  by induction on the
dimension of $N$. This is trivial when $N$ is one-dimensional. Now
consider the general case, assuming the validity in lower-dimensional
cases. If the $\Phi$-action on $N/\Gamma$ is minimal then by
Proposition~\ref{minimal} it follows that $N$ is abelian, as desired.
Now suppose that the $\Phi$-action is not minimal.  Let
$V=N/[N,N]$ and  $\Delta =\Ga [N,N]/[N,N]$. Since the $\Phi$-action is
not minimal, by Parry's theorem the factor action on $V/\Delta$
is also not minimal. This implies that there exists a proper vector
subspace
$W$ of $V$ containing $\Phi [N,N]$ such that $W\Delta$ is closed.  Let
$M$ be the Lie subgroup of $N$ containing $[N,N]$ such that
$W=M/[N,N]$. Then $\Phi$ is contained in $M$ and the $\Phi$-action on
$M\Ga/\Ga$ is almost automorphic. By the induction hypothesis there
exists a closed connected  abelian subgroup $A$ containing  $\Phi$
such that $A\Ga$ is closed, thus proving (iii).

iii) $\implies$ i): When (iii) holds $A\Ga/\Ga$ is a closed invariant
subset the restriction to which is equivalent to a
translation on $A/A\cap \Ga$ which is a torus, and in particular it is
an isometry with respect to a metric. Hence $\Ga$ is almost
automorphic.  This completes
the proof. \hfill $\Box$.

\section{Proofs  of the main results}\label{mainthm}

In this section we shall deduce from the case of the translations the
general case of the main theorem. The main point is to realise the
unipotent affine automorphisms as restrictions of translations of
higher dimensional compact nilmanifolds. This corresponds to
considering the ``suspension'' of the given affine automorphism, in
the sense of dynamics. We begin by recalling the construction.

Let $N$ be a simply connected nilpotent Lie group, $\Ga$  be a lattice
in $N$, $U$ a unipotent automorphism of $\N/\Ga$ and $\tilde U$ the
corresponding automorphism of $N$. Then there exists a one-parameter
subgroup $\Phi=\{\tilde U_t\}$ consisting of unipotent automorphisms
of $N$, with $\tilde U=\tilde U_1$. Let $M=\Phi \cdot N$, the
semidirect product of $\Phi$ and $N$ with respect to the action of
$\Phi$ by automorphisms.  Then $M$ is a simply connected nilpotent Lie
group. Since $\tilde U (\Ga)=\Ga$, it follows that $\Delta:=\{\tilde
U^n\mid n\in \Z\}\cdot \Ga$ is a subgroup of $M$. Furthermore $\Delta$
is a lattice in $M$. We realise $N/\Ga$ canonically as a subset of
$M/\Delta$. We note also that the $U$-action on $N/\Ga$ is the
restriction of the translation of $M/\Delta$ by the element $\tilde
U$; for $g\in \N$, $(\tilde U\cdot g)\Delta=\tilde U(g)\tilde U
\Delta= \tilde U(g)\Delta$ and the latter is the same as $U(g)\Ga$
under the realisation of $N/\Ga$ as a subset of $M/\Delta$.

\medskip
\noindent{\it Proof of Theorem~\ref{thm:main}}: We follow the
  notation as in the hypothesis and as introduced above. Suppose that
  $T$ is almost automorphic. We have
  realised $N/\Ga$ as a subset of $M/\Delta$, such that the $U$-action
  on $N/\Ga$ is the restriction of a translation of $M/\Delta$, by an
  element, specifically $\tilde U$ viewed as an element of $M$. Then
  $T$ is also the restriction of a translation, by the element
  $a\tilde U \in M$, and we note that $\Ga$ being almost automorphic for
  the $T$-action implies that $\Delta$ is almost automorphic for the
  translation by $a\tilde U$. Hence there exists a closed connected
  abelian subgroup $B$ of $M$, containing $a\tilde U$, such that
  $B\Delta $ is   closed. The image of $B$ modulo $N$ is a connected
  subgroup of $\Phi$ containing $\tilde U$ and hence it is the whole
  of $\Phi$. Therefore $B=\Phi \cdot A$ where $A=B\cap N$. Also since
  $B\Delta$ is closed
and $\Phi \cap \Delta$ is a  lattice in $\Phi$ it follows
that $A\Ga$ is closed.  Since $a\tilde U\in B$ we get $a\in A$. Also
any $g\in A\subset B$ commutes with $\tilde U$ as an element of $M$,
which shows that $\tilde U(g)= \tilde U g \tilde U^{-1}=g$. This
proves the ``only if'' part of the theorem. The ``if'' part is obvious
since under the condition the $T$-orbit of $\Ga$ is contained in
$A\Ga/\Ga$ which is a torus, and the action is a translation of the
torus and hence an isometry with respect to a metric. This proves the
theorem. \hfill $\Box$

\bigskip \noindent
{\it Proof of Corollary~\ref{cor:main}} :
Suppose that $T$ is almost automorphic. Let  $x\in N$ be
arbitrary. Then  $x\Ga$
is almost automorphic for the $T$-action.  We view $N/\Ga$ as $N/x\Ga
x^{-1}$, choosing $x\Ga$ as the base point. Then the $T$-action
corresponds to $T_{a_x}\circ U_x$, where $a_x=a\tilde U(x)x^{-1}$, and
$U_x$ is the automorphism of $N/x\Ga x^{-1}$ induced by the
automorphism $\tilde U_x$ given by $\tilde U_x(g)= x\tilde
U(x)^{-1}\tilde U(g)\tilde U(x)x^{-1}$. Applying
Theorem~\ref{thm:main}  we get that there
  exists a closed connected subgroup $A_x$ such that $A_xx\Ga$ is
  closed, $a_x\in A_x$ and $\tilde U_x(g)=g$ for all $g\in A_x$.

Let $\cal A$ be the set of all closed connected abelian subgroups $A$
such
that $A\Ga$ is closed. Then $\cal A$ is countable; this follows from
the fact that the Lie subalgebra corresponding to $A\in \cal A$ is
determined by linear equations with rational coefficients. We have
$x^{-1}A_xx\in \cal A$ for all $x\in N$ and hence it follows that
there exists $A\in \cal A$ such that the set $E$ consisting of all $x$
such that $x^{-1}A_xx= A$ (which is  closed and hence a Borel subset),
has positive Haar measure. We have $a_x\in A_x$ and
hence $x^{-1}a\tilde U(x)\in x^{-1}A_xx= A$ for all $x\in E$. Since the set
of $x$ such that  $x^{-1}a\tilde U(x)\in  A$ is a submanifold of $N$ and $E$
has positive Haar measure, this
implies that $ x^{-1}a\tilde U(x)\in  A$ for all $x\in N$.

Now consider any
$x\in E$. We have $\tilde U_x(g)= x\tilde
U(x)^{-1}\tilde U(g)\tilde U(x)x^{-1}=g$ for all $g\in A_x$. Thus
$\tilde U(x^{-1}gx)=x^{-1}gx$ for all $g\in A_x$, and since
$A=x^{-1}A_xx$ this means that $\tilde U(g)=g$ for all $g\in A$. This
shows that $A$ has the desired properties and thus proves the ``only
if'' part of the corollary.

Now suppose that there exists a subgroup $A$ as in the statement. Let
$x\in N$ be arbitrary. Then $xA\Ga/\Ga$ is a closed subset of
$N/\Ga$. For any $g\in A$ we have $T(xg\Ga)=a\tilde U(xg)\Ga=a\tilde
U(x)\tilde U(g)\Ga=a\tilde U(x)g\Ga$, since $\tilde U(g)=g$, and
$a\tilde U(x)\in xA$ by hypothesis, so we get that $T(xg\Ga)\in
xA\Ga/\Ga$. Thus $xA\Ga/\Ga$ is $T$-invariant. Furthermore the above
relation also shows that under the canonical correspondence of $xA\Ga/\Ga$
with $A\Ga/\Ga$ the restriction of $T$ to the former corresponds to
the translation action by the element $x^{-1}a\tilde U(x)$. Now
$A\Ga/\Ga$ is a torus and the translation is an isometry with respect
to a metric on it, so the map is almost automorphic. Since this holds
for every $x$ it follows that $T$ is almost automorphic.

Finally, we note that if $T$ is also minimal then $A\Ga/\Ga$ as above
has to be the whole of $N$, which implies that $N$ is abelian, and
that $T$ is a translation of $N/\Ga$. \hfill $\Box$

\bigskip
\noindent{\it Proof of Corollary~\ref{transl}} : The ``if'' part is
straightforward and we omit the proof. Now suppose $T_a$ is automorphic.
 Then by Corollary~\ref{cor:main} there exists a closed connected
 abelian subgroup $A'$ such that $A'\Ga$ is closed and $x^{-1}ax\in A'$
  for all $x\in N$. Let $B$ be the smallest closed subgroup of $N$
  containing $\{x^{-1}ax \mid x\in N\}$. Then $B$ is a closed normal
  subgroup of $N$ contained in $A'$. Hence $\overline {B\Ga}$ is a closed
  subgroup of $A'\Ga$. Let $A$ be the connected component of the identity
  in $\overline {B\Ga}$. Then $A$ is a closed connected subgroup.
  Moreover, as $\overline{B\Ga}\subset A'\Ga$ and the latter is a closed
  subgroup, $A\subset A'$ and hence $A$ is abelian. Since
  $A\Ga$ is the closed subgroup $\overline {B\Ga}$ and $A$ is identity
  component, $A$ is normalised by $\Ga$. Since $\Ga$ is a lattice in a
  simply connected nilpotent Lie group this implies that $A$ is normal in
  $N$ (see \cite{raghunathan}, Corollary 2 of Theorem 2.3). This completes the proof. \hfill $\Box$

\bigskip
\noindent{\it Proof of Corollary~\ref{cor:Lie}} :  By
Corollary~\ref{cor:main}, for every $x\in N$, $x^{-1}a\tilde
U(x)$ and $a$ are  fixed by $\tilde U$, and hence so is $x^{-1}a\tilde
U(x)a^{-1}$. Given $\xi \in \mathfrak N$, applying this to $\exp t\xi$,
$t\in \R$, in place of $x$, and
differentiating with respect to $t$ we get
that  $\Ad_a\circ \mathfrak U (\xi)-\xi$ is fixed by $\mathfrak U$. Thus
$\mathfrak U (\Ad_a\circ \mathfrak U (\xi)-\xi)=
\Ad_a\circ \mathfrak U (\xi)-\xi$ for all $\xi \in \mathfrak N$, that is,
$(\mathfrak U -I)(\Ad_a\circ \mathfrak U -I)=0$. This proves (i)

Also, by
Corollary~\ref{cor:main}, the elements $x^{-1}a\tilde
U(x)a^{-1}$ commute with each other. Let $\mathfrak M=(\Ad_a\circ \mathfrak U
-I)(\mathfrak N)$. Then $\mathfrak M$ is readily seen to be a Lie subalgebra
of $\mathfrak N$. Given $\xi \in \mathfrak M$ it is
tangent to the curve  $x_t^{-1}a\tilde
U(x_t)a^{-1}$ where $x_t=\exp t\xi$ for all $t\in \R$. Since for $\xi,
\eta \in \mathfrak M$ the elements in the corresponding curves as above
commute with each other, it follows that $\xi$ and $\eta$ commute in
$\mathfrak N$. Hence $\mathfrak M$ is commutative. This proves the
Corollary. \hfill $\Box$

\section{Some further consequences}

In this section we discuss some auxiliary results around the theme
of the main results.

For $N=\R^d$, $d\geq 1$ we have the following characterisation of
almost automorphic affine automorphisms, which may be compared with
Corollary~\ref{cor:Lie}.

\begin{proposition}\label{prop:R^d}
Let $N=\R^d$, $d\geq 1$, and $\Ga= \Z^d$. Let $T=T_a\circ U$ be
an affine automorphism of $N/\Ga$, where $a\in N$ and $U$ is an
automorphism of $N/\Ga$ induced by a unipotent automorphism $\tilde U$
of $N$. Then $T$ is almost automorphic if and only if
$(\tilde U-I)^2=0 $ (here $I$ and $0$ are the identity and zero
matrices respectively)  and $\tilde U(a)=a$.
\end{proposition}

\proof Suppose $T$ is almost automorphic. Then by
Corollaries~\ref{cor:Lie} and \ref{cor:main} $(\tilde U-I)^2=0 $ and
$\tilde U(a)=a$. Conversely suppose $(\tilde U-I)^2=0 $
 and $\tilde U(a)=a$. Then  every $T$-orbit is contained in
an orbit, say $O$, of the subgroup $\{x \in \T^d \mid U(x)=x\}$ and
restricted to $O$ the $U$-action is equivalent to a translation.
Hence $T$ is almost automorphic. \hfill $\Box$

\medskip
The following example shows that the converse of
Corollary~\ref{cor:Lie} is however not true in
general, even for automorphisms; (for the affine case we already have
another necessary condition, that $a$ is fixed by $\tilde U$, not
incorporated in Corollary~\ref{cor:Lie} focussing on the
 conditions at the Lie algebra level.

\begin{example}{\rm Let $\mathfrak N$ be the $4$-dimensional Lie algebra
with a set of basis vectors $\{\xi_1,\xi_2, \xi_3, \xi_4\}$ and the
Lie bracket
satisfying the conditions $[\xi_1, \xi_2]=\xi_3=-[\xi_2, \xi_1]$,
$[\xi_2,\xi_3]=\xi_4=-[\xi_3, \xi_2]$ and $[\xi_i, \xi_j]=0$ for $(i,j)$
different from $(1,2), (2,1), (2,3)$ and $(3,2)$. It can be seen that
$\mathfrak N$ is a nilpotent Lie algebra. Let $\mathfrak U: \mathfrak N\to \mathfrak
N$ be the linear map such that $\mathfrak U (\xi_1)=\xi_1+\xi_2$ and
$\mathfrak U (\xi_i)=\xi_i$ for $i=2,3$ and $4$. It can be verified that
$\mathfrak U$ is a Lie algebra  automorphism and $(\mathfrak U -I)^2=0$.
Also the image of $\mathfrak U -I$ is the (one-dimensional) span of
$\xi_2$, which is an abelian Lie subalgebra.
Let $N$ be the simply connected Lie group corresponding
to $\mathfrak N$ and let $\tilde U$ be the Lie automorphism of $N$
corresponding to $\mathfrak U$. Let $\Lambda$ be the subgroup of $\mathfrak N$
generated by the basis  $\{\xi_1,\xi_2, \xi_3, \xi_4\}$ and $\Ga =\exp
(\Lambda)$,  where
$\exp$ denotes the exponential map of $\mathfrak N$ onto $N$. Since the
basis has rational structure constants it follows that $\Ga$
is a lattice in $N$ (cf. \cite{raghunathan}, Theorem 2.12). Also, since
$\Lambda $ is $\mathfrak U$-invariant, it follows that $\Ga$ is $\tilde
U$-invariant. Let $U$ be the corresponding automorphism of
$N/\Gamma$. Let $M$ be the
connected subgroup of $N$ corresponding to the Lie subalgebra spanned
by $\{\xi_2,\xi_3, \xi_4\}$; then $M$ is also the set of fixed points
of $\tilde U$. For all $t\in \R$ let $x_t=\exp t\xi_1$, and consider
the subset $E=\{x_t^{-1}\tilde U(x_t)\mid t\in \R\}$ of $M$. We see
that $E$ generates $M$; this may be proved by observing that in
$M/[M,M]$, which may be viewed canonically as $\R^2$, the image of $E$
is an affine line. In particular $E$ is not contained in an abelian
subgroup. Hence by Corollary~\ref{cor:main} we get that
$U$ is not almost automorphic.~\hfill $\Box$}
\end{example}

\medskip
The phenomenon involved in the above Example is seen in a more general
form in the following proposition.

\begin{proposition}\label{proposition}
Let $\mathfrak N$ be a nilpotent Lie algebra generated, as a Lie algebra,
by two elements $\xi$ and $\eta$, and $\mathfrak U$ be a unipotent Lie
automorphism of $\mathfrak N$ such that  $\mathfrak U (\xi)=\xi+\eta$.
Let $N$ be a simply connected nilpotent Lie group with  $\mathfrak N$
as its Lie algebra. Let $\tilde U$ be
the Lie automorphism of $N$ corresponding to $\mathfrak U$. Let  $\Ga$ be
a  lattice in $N$ invariant under $\tilde U$   and let $U$ be
the automorphism of $N/\Ga$ induced by $\tilde U$. If $U$ is almost
automorphic then $N$ contains a $(d-1)$-dimensional abelian $\tilde
U$-invariant normal subgroup $M$ such that $N$ is the semidirect
product of $\{\exp t\xi\}$ and $M$, $M\cap \Ga$ is a lattice in
$M$, and the $\tilde U$-action on $M$ is trivial.
\end{proposition}

\proof  Since $\mathfrak N$ is generated by $\xi $ and $\eta$
it follows that $\mathfrak N/[\mathfrak N, \mathfrak N]$ is a two dimensional vector
space spanned by the images of $\xi$ and $\eta$. Since $\mathfrak U
(\xi)=\xi+\eta$ this implies that $\langle \eta \rangle +[\mathfrak N,
\mathfrak N]$ is a $\mathfrak U$-invariant Lie ideal in $\mathfrak N$. Let $\mathfrak
M=\langle \eta \rangle +[\mathfrak N, \mathfrak N]$ and $M$ be the simply
connected subgroup of $N$ corresponding to $\mathfrak M$. Then $M$
is a $(d-1)$-dimensional $\tilde U$-invariant subgroup  and  $N$ is
the semidirect product of $\{\exp t\xi\}$ and $M$. Also, since
$[N,N]\cap \Ga$ is a lattice in $[N,N]$ and $M/[N,N]$ is the set of
fixed points of  the factor of $\tilde U$ on $N/[N,N]$ it follows that
$M\cap \Ga$ is a lattice in $M$. We conclude the
proof by showing that $M$ is abelian and the $\tilde U$ action on $M$
is trivial.

For all $t\in \R$ let $x_t=\exp t\xi$. Then as seen before, for every
$t\in \R$, $x_t M\Gamma/\Gamma$ is $U$-invariant and when it is viewed as
$M/(M\cap\Gamma)$ the $U$-action on
$x_t M \Gamma/ \Gamma$, $t \in \mathbb{R}$ corresponds to the
affine automorphism induced by $T_{a_t} \circ \tilde U$, where
$a_t={x_t}^{-1}\tilde U(x_t)=
\exp(-t\xi) \exp(t(\xi+\eta))$. Let $L=M/[M,M]$ and
$\Delta =[M,M](M\cap \Ga)$. Let $\mathfrak L$ be the Lie algebra of $L$
and $\tau :\mathfrak L \to \mathfrak L$ be the map defined by
$\tau (\mu +[\mathfrak M,\mathfrak M])=[\xi, \mu]+[\mathfrak M,\mathfrak M]$ for all
$\mu\in \mathfrak M$. Since $\mathfrak N$ is generated by $\xi$ and $\eta$ it
follows that $\mathfrak L$ is  spanned by $\{\eta, \tau (\eta), \dots,
\tau^{n-1}(\eta)\}$, where $n$ is the least positive integer such that
$\tau^n(\eta)=0$, and furthermore the set forms a basis of $\mathfrak L$.

Now $N/[M,M]$ may be viewed as the semidirect product of the
one-parameter subgroup $\Phi:=\{x_t\}$ and $\mathfrak L$, identifying $L$
with $\mathfrak L$ canonically. In turn the group
may be realised as a group of  $(n+1)\times (n+1)$ matrices, with
$x_t$ represented by $\left(\begin{matrix} \exp t\nu&0\\ 0&1\end{matrix}\right)$
%$\left(
%      \begin{array}{cc}
%      \exp t\nu &  0\\
%      0& 1\\
%     \end{array} \right)$,
where $\nu=(\nu_{ij})$ is the $n\times n$ matrix given by $\nu_{ij}=1$
if $j=i+1$ and $0$ otherwise, and the elements of $\mathfrak L$ are
represented by $\left(\begin{matrix} I_n&v\\ 0&1\end{matrix}\right)$
%$\left(
%      \begin{array}{cc}
%      I_n &  v\\
%      0& 1\\
%      \end{array} \right)$,
with $I_n$ the $n\times n$ identity matrix and $v$ a column vector
$^t(v_1, \dots, v_n)$; specifically the
entries $v_1, \dots, v_n$ are given by the coordinates of the element
with respect to the basis $\{\tau^{n-1}(\eta), \dots, \eta\}$ (in that
order).
From this the image of $a_t$ in $L$ can be computed explicitly to be
$\sum \limits_{k=0}^{n-1} c_k t^{k+1} \tau^k(\eta)$,
where $c_k={(-1)^{k}}/{k!}$ for all $k$.

We can find $t$ such that $1, c_0 t, c_1 t^2,\ldots c_{n-1} t^n$ are
linearly independent over $\mathbb{Q}$. This is because all $c_k$
being nonzero, the set of $t$ for which the elements are linearly
dependent is countable. For $t$ satisfying this condition, as
discussed above, the action
$T_{a_t}\circ U$ on $L/\Delta$ is minimal, and in turn
the $T_{a_t}\circ U$ action on $M/(M\cap\Gamma)$ is minimal. Since it
is almost automorphic by
Corollary~\ref{cor:main} that $M$ is abelian and the $\tilde U$ action on
$M$ is trivial.
 \hfill $\Box$

\section{Miscellania}

In this section we discuss the issue in a broader context. Firstly we
prove the following result which extends the scope of applicability of the
results of the earlier sections.

We shall say that a homeomorphism $\varphi $ of a locally compact
second countable
space $X$ {\it  admits a convergent trajectory} if there exists $x\in
X$ such that $\varphi (x)\neq x$ and either $\{\varphi^k(x)\}_{k=0}^\infty$ or
$\{\varphi^{-k}(x)\}_{k=0}^\infty$ converges in $X$, as $k\to
\infty$; we note that in either case the limit is a fixed point
of $\varphi$. Clearly, if $\varphi$ admits a convergent trajectory, or
more generally if a factor of $\varphi$ admits a convergent
trajectory, then $\varphi$ is not almost automorphic. This means in
particular that if $T$ as in the
hypothesis of the next proposition is also assumed to be  almost
automorphic, then
only the second alternative as in the conclusion can hold.

\begin{proposition}\label{general}
Let $N$ be a simply connected nilpotent Lie group and $\Gamma$
be a lattice in $N$.   Let $T=T_a\circ  A$ be an affine
automorphism of  $N/\Ga$, where $a\in N$ and $A$ is an automorphism of
$N/\Ga$. Let $\tilde A$ be the automorphism corresponding to $A$.
Then one of the following holds:

i) there exists a closed connected normal $\tilde A$-invariant
subgroup $M$ of $N$ such that $M\Ga$ is closed, and the factor
$\bar T$ of $T$ on $N/M\Ga$ admits a convergent trajectory.

ii)  there exists $r\in \N$ such that $\tilde A^r$ is unipotent.
\end{proposition}

\proof
Firstly consider the case when $N=\R^d$, for some $d\geq 1$.
Let $V$ be the largest $\tilde A$-invariant subspace of $\R^d$ on which
 $1$ is the only eigenvalue, and $W=\R^d/V$. Since $\tilde A$ leaves
 $\Ga$ invariant,
 it follows that the image $\Delta=(V+\Ga)/\Ga$ of $\Ga$ in $W$ is a lattice
 in $W$. It suffices to prove the proposition for the factor of
$T$ on $R/W$ and hence without loss of generality we may
 assume that $V$
 is trivial and $W=\R^d$; that is, $1$ is not an eigenvalue of $\tilde
 A$. The last condition implies,  via simple linear algebra, that $T_a\circ
 \tilde A$ has a fixed point $b\in  \R^d$. Then
 $T_a\circ  \tilde A= T_b\circ\tilde A\circ T_b^{-1}$.
 If $\tilde A$ has an eigenvalue of absolute value different from $1$
 then $\tilde A$ admits a convergent trajectory and hence so does  $T_a\circ
 \tilde A$, and so statement~(i) holds in this case. Now suppose that
 all the eigenvalues of $\tilde A$ are of absolute value $1$. Since $\tilde A$
 leaves $\Ga$ invariant this further implies that all the eigenvalues
 are roots of unity;  (one shows
this by considering the characteristic polynomials of $\tilde A^k$, $k\in
\N$, and showing that under the given condition their roots, which are
eigenvalues of
$\tilde A^k$, would belong to a fixed finite set). Hence there exists
$r\in \N$ such that $\tilde A^r$ is unipotent. This proves the
proposition when $N$ is abelian.

Now consider the general case. Let $\bar A$ denote the factor of $\tilde A$
on $N/[N,N]$. Suppose that alternative (i) does not hold. From the special
case we then get that there
exists $r\in \N$ such that $\bar A^r$ is unipotent. Therefore in
proving the proposition we may without loss of generality assume that
$\bar A$ is unipotent. By an argument well-known to Lie group theorists
this implies that  $\tilde A$
is also unipotent. Let us recall the argument here, briefly, for the reader's
convenience:  Let $\mathfrak N$ be the Lie algebra of $N$ and $\mathfrak
M$ be the largest $d \tilde A$-invariant subspace on which $d\tilde A$
is unipotent. It can be seen that $\mathfrak M$ is a Lie subalgebra of
$\mathfrak N$, and as $\bar A$ is unipotent  we get $\mathfrak M +[\mathfrak N,
\mathfrak N]=\mathfrak N$. Substituting from the equation successively we get
$\mathfrak N=\mathfrak M+ [\mathfrak M +[\mathfrak N, \mathfrak N], \mathfrak M +[\mathfrak N,\mathfrak
N]]\subseteq  \mathfrak M + [\mathfrak N, [\mathfrak N, \mathfrak N]] \subseteq \cdots
\subseteq \mathfrak M$, since $\mathfrak N$ is nilpotent. Thus $\mathfrak N=\mathfrak M$ which
shows that $\tilde A$ is unipotent. This completes the proof of the
proposition. \hfill $\Box$

\medskip
\begin{remark}
{\rm
We recall that a homeomorphism $\varphi$ of a locally compact second
countable space $X$ is said to be {\it distal} if for  $x,y\in X$,
$x\neq y$, there does not
exist any sequence $\{k_i\}$ in $\Z$ such that $\varphi^{k_i}(x)$ and
$\varphi^{k_i}(y)$ converge to the same
point as $i\to \infty$. Proposition~\ref{general} also shows that if
an affine
automorphism $T=T_a\circ A$, as in the statement of the proposition, is
assumed to be distal then $A^r$ is unipotent for some $r\in \N$; in
other words,
the automorphism part of $T^r$ is unipotent. Related general
results for distal group actions by group automorphisms, with analogous
conclusions about the action of a subgroup of finite index, may be found in
\cite{abels}. The preceding conclusion can indeed be deduced, albeit not
directly, from the results of \cite{abels}. On the other hand
Proposition~\ref{general} provides a simple and direct argument
in the case at hand.
 }
\end{remark}

We conclude with the following remark.

\begin{remark}
{\rm The notion of almost automorphic homeomorphisms extends naturally
  to group actions $G\times X\to X$, where $G$ is a group acting by
  homeomorphisms of $X$: the action is almost automorphic if for any
sequence $\{g_i\}$ in $G$ and $x, y\in X$ if $g_ix\to y$ then
$g_i^{-1}y\to x$. Clearly when an action is almost automorphic then
the homeomorphism corresponding to the action of any individual element
is almost automorphic. It seems to us however that the converse of
this is not true. This would follow if it is proved, which also we
expect to be true, that any group of automorphism of a
compact nilmanifold whose action is almost automorphic is necessarily
{\it abelian}; (it can be readily seen that there exist nonabelian
subgroups of $SL(4, \Z)$ all whose elements are unipotent elements of nilrank
at most $1$, so the action of each of them on $\T^4$ is almost
automorphic (by Proposition~\ref{prop:R^d})).

}

\end{remark}
\noindent
{\it Acknowledgement}: The first and third named author would like to thank
the School of Physical Sciences (SPS), Jawaharlal Nehru University (JNU),
New Delhi, India for hospitality while some of the work was done. The third
named author would also like to thank the National Board of Higher Mathematics
(NBHM), DAE, Govt.\ of India for a post-doctoral fellowship. The second named
author would like to thank the School of Mathematics, Tata Institute of
Fundamental Research (TIFR), Mumbai, India for hospitality while some of
the work was done.
%\end{section}

\bigskip

\begin{flushleft}
S.G.~Dani\\
Department of Mathematics\\
Indian Institute of Technology Bombay\\
Powai, Mumbai 400 076\\
India\\
{\tt shrigodani@gmail.com}\\[.7cm]

Riddhi Shah\\
School of Physical Sciences\\
Jawaharlal Nehru University\\
New Delhi 110 067\\
India\\
{\tt riddhi.kausti@gmail.com}\\[.7cm]

Puneet Sharma\\
COE - Systems Science\\
Indian Institute of Technology Jodhpur\\
Residency Road, Ratanada\\
Jodhpur - 342 011\\
India\\
{\tt puneet8111.iitd@gmail.com}

\end{flushleft}

\end{document}